\documentclass[10pt,english]{amsart}
\usepackage{amsmath,amssymb,amsthm,amsfonts,graphicx,color}
\usepackage{amssymb,geometry}

\usepackage[colorlinks=true]{hyperref}

\makeatletter
\def\@currentlabel{2.1}\label{e:dispaa}
\def\@currentlabel{2.21}\label{e:dispau}
\def\@currentlabel{2.22}\label{e:dispav}
\def\@currentlabel{2.23}\label{e:dispaw}
\def\@currentlabel{2.24}\label{e:dispax}
\def\theequation{\thesection.\@arabic\c@equation}
\makeatother

\hypersetup{linkcolor=black,citecolor=black,filecolor=black,urlcolor=black} 

\definecolor{dullmagenta}{rgb}{0.4,0,0.4}   
\definecolor{darkblue}{rgb}{0,0,0.4}

\newcommand{\R} {\mathbb R}
\newcommand{\cuad}{{\sqcap\kern-.68em\sqcup}}

\newcommand{\be}{\begin{equation}}
\newcommand{\ee}{\end{equation}}

\newcommand{\sech}{\mathop{\mbox{\normalfont sech}}\nolimits}

\makeatletter
\newcommand*\rel@kern[1]{\kern#1\dimexpr\macc@kerna}
\newcommand*\widebar[1]{%
  \begingroup
  \def\mathaccent##1##2{%
    \rel@kern{0.8}%
    \overline{\rel@kern{-0.8}\macc@nucleus\rel@kern{0.2}}%
    \rel@kern{-0.2}%
  }%
  \macc@depth\@ne
  \let\math@bgroup\@empty \let\math@egroup\macc@set@skewchar
  \mathsurround\z@ \frozen@everymath{\mathgroup\macc@group\relax}%
  \macc@set@skewchar\relax
  \let\mathaccentV\macc@nested@a
  \macc@nested@a\relax111{#1}%
  \endgroup
}
\makeatother

\renewcommand{\theequation}{\thesection.\arabic{equation}}
 
 \newtheorem{lemma}{Lemma}[section]
\newtheorem{definition}{Definition}
\newtheorem{theorem}{Theorem}[section]

\newtheorem{remark}{Remark}[section]
\newcommand{\bremark}{\begin{remark} \em}
\newcommand{\eremark}{\end{remark} }

\begin{document}

\title[Free boundary problems and the  Liouville equation]
{Free boundary problems arising in  the theory of maximal solutions of equations with exponential nonlinearities}
%
%


\author{Micha{\l } Kowalczyk}
\address{Departamento de Ingenier\'{\i}a Matem\'atica and Centro
de Modelamiento Matem\'atico (UMI 2807 CNRS), Universidad de Chile, Casilla
170 Correo 3, Santiago, Chile.}
\email {kowalczy@dim.uchile.cl}

\author{Angela Pistoia}\address{Dipartimento SBAI, Sapienza Universit\'a di Roma, via Antonio Scarpa 16, 00161, Roma, Italy.}
\email {angela.pistoia@uniroma1.it}

\author{Piotr Rybka}
\address{Institute of Applied Mathematics and Mechanics, The University of Warsaw, Banacha 2, 02-097 Warsaw, Poland.} \email{rybka@mimuw.edu.pl}

\author{Giusi Vaira}\address{Dipartimento di Matematica e Fisica,   Universit\'a della Campania ``L. Vanvitelli", viale Lincoln 5, 81100, Caserta, Italy.}
\email {giusi.vaira@unicampania.it}

\thanks{M. Kowalczyk was partially supported by Chilean research grants Fondecyt 1130126 and 1170164 and Fondo Basal AFB170001 CMM-Chile. Part of this work was done during his visits at the University of Warsaw and Hiroshima University, they were supported by the program "Guests" of the Warsaw Center for Mathematical Sciences. A. Pistoia was partially supported by Sapienza   research grant  ``Nonlinear PDE's in geometry and physics''. A part of P.Rybka's work was done during his visit to the Universidad  de Chile, whose hospitality is gratefully acknowledged. G. Vaira was partially supported by GNAMPA research grant ``Esistenza e molteplicit\'a di soluzioni per alcuni problemi ellittici non lineari".}
\begin{abstract}
We  consider equations of the form  $\Delta u +\lambda^2 V(x)e^{\,u}=\rho$ in various  two dimensional settings.  We assume that $V>0$ is a given function, $\lambda>0$ is a small parameter and $\rho=\mathcal O(1)$ or $\rho\to +\infty$ as $\lambda\to 0$. In a recent paper \cite{kpv2018}  we proved the existence of the maximal solutions for a particular choice $V\equiv 1$, $\rho=0$ when the problem is posed in doubly connected domains  under Dirichlet boundary conditions. We related the maximal solutions with a novel free boundary problem. The purpose of this note is to derive the corresponding free boundary problems in other settings. Solvability of such problems is, viewed formally, the necessary condition for the existence of the maximal solution.    
\end{abstract}
\date{}\maketitle


\section{The Liouville equation }
\setcounter{equation}{0}
\subsection{The maximal solution}

The following problem is known as  the Liouville equation:
\begin{equation}
\label{liu1} 
\begin{aligned}
\Delta u+\lambda^2 e^{\,u}&=0, \quad\mbox{in}\ \Omega\subset \R^2,\\
u&=0, \quad\mbox{on}\, \partial\Omega,
\end{aligned}
\end{equation}
In \cite{nagasaki1990} Nagasaki and Suzuki proved that for any sequence  of 
$\lambda_n\to 0$  one of the following  holds  for the sequence of the  corresponding solutions $u_n$ of (\ref{liu1}):
\begin{itemize}
\item[(i)] $\|u_{n}\|_{L^\infty(\Omega)}\to 0$.
\item[(ii)] Solutions $u_{n}$ blow up at the set of isolated points in $\Omega$.
\item[(iii)] Blow up occurs in the whole $\Omega$ i.e. $u_{n}(x)\to \infty$ for all $x$.
\end{itemize}
In the case (i) we speak of the minimal solutions.  They can be obtained by minimizing the functional
\[
\int_\Omega\frac{1}{2} |\nabla u|^2-\lambda^2 \int_\Omega e^u
\]
which is coercive, due to the Trudinger-Moser inequality, 
when $\lambda$ is small enough. In the case (ii) we speak of the bubbling solutions. They  blow up at $k$ isolated points $a_j\in \Omega$. 
Near each $a_j$, $j=1, \ldots, k$,  the local profile of the bubble is a scaling of 
\begin{equation}
\label{stndr bubl}
w(r)=\log\frac{8}{(1+r^2)^2}
\end{equation}
which is a solution of $\Delta w+e^{\,w}=0$ in $\R^2$, we call it  the standard bubble. 
When $x\approx a_j$ we have  $u(x)\approx w(|x-a_j|/\lambda)-4\log\lambda$,   so  that 
\begin{equation}
\label{mass1} 
\lambda^2 \int_\Omega e^{\,u}\to 8\pi k, \quad \mbox{as}\ \lambda \to 0,
\end{equation}
while for the minimal solution this last limit is $0$. Finally, when the last alternative happens we speak of the {maximal solutions}. Because of the well known result of Br\'ezis and Merle   \cite{doi:10.1080/03605309108820797} (see also \cite{MR1322618}) we know that for (iii) to hold we must have
\[
\lim_{m\to \infty}\lambda^2_{n} \int_\Omega e^{\,u_{n}}=\infty.
\]
From \cite{NAGASAKI1990144} we know that at least in the case of the annulus all three alternatives may occur. In particular solutions satisfying (i) or (iii) are radial and those satisfying (ii) have $k$ fold symmetry. Constructions of the bubbling solutions can be found in \cite{Baraket1997, MR2157850} (for a similar result for the mean field equation see \cite{Esposito2005}). Summarizing these results  we know that for any multiply connected domain there exists a solution with $k$ bubbles, $k\in \mathbb N$.

From now on we  suppose  that  $\lambda>0$ is a small parameter and $\Omega$ is a bounded, smooth domain, which is doubly connected and such that the bounded component of $\R^2\setminus \Omega$ is not a point.
To state the existence result proven in \cite{kpv2018} we will introduce the notion of the harmonic measure of a closed curve $\gamma\in\Omega$.
\begin{definition}\label{liu def 1}
Let $\gamma$ be a smooth, simple closed curve in $\Omega$. A function $H_\gamma\in {C}^{2}(\Omega\setminus \gamma)\cap {C}^0(\Omega) $ is called the harmonic measure of $\gamma$ in $\Omega$ if the following holds:
\begin{equation}
\begin{aligned}
\label{liu 3}
\Delta H_\gamma&=0,\quad \mbox{in}\ \Omega\setminus \gamma,\\
H_\gamma&=0, \quad \mbox{on}\ \partial\Omega,\\
H_\gamma&=1, \quad\mbox{on}\ \gamma.
\end{aligned}
\end{equation}
\end{definition}
Since $\gamma$ is a simple closed curve it divides the plane, and consequently $\Omega$, into two disjoint components, 
{$\Omega\setminus\gamma = \Omega^+\cup \Omega^-$. By definition, $\Omega^+$ is contained in the bounded component of $\R^2\setminus\gamma$,
it can be called interior with respect to $\gamma$. 
The other component, $\Omega^-$, can be called the exterior with respect to $\gamma$.}  We will also set:
\begin{align*}
H^\pm_\gamma=H_\gamma\left|_{\Omega^\pm}\right. .
\end{align*}
Functions $H^\pm_\gamma$ are harmonic in their respective domains $\Omega^\pm$ and they satisfy homogeneous Dirichlet boundary conditions on $\partial\Omega^\pm\setminus\gamma$. Finally, for future purpose by $n$ we will denote the unit normal vector on $\gamma$. 
{It defines the orientation so that} 
$n$ is the exterior  unit normal of $\Omega^-$ on $\gamma$, and  the interior unit normal of $\Omega^+$ on $\gamma$. 
{Keeping this in mind,} we have \cite{kpv2018}:
\begin{lemma}\label{lemma free}
Let $\Omega\in \R^2$ be a bounded,  smooth, doubly connected set  such that the bounded component of $\R^2\setminus \Omega$  is not a point. There exists a simple, closed and smooth curve $\gamma\in \Omega$ such that its harmonic measure satisfies
\begin{align}
\label{liu 4}
\partial_n H_\gamma^++\partial_n H^-_\gamma=0 \quad \mbox{on}\ \gamma.
\end{align}
\end{lemma}

By the generalization of the Riemann mapping theorem for multiply connected domains there exists a  holomorphic, bijective map $\psi\colon\Omega\to B_{R_1}\setminus B_{R_2}$ with some $R_1>R_2>0$. {We will show in Section 1.2} that $\psi(\gamma)=\partial B_{R}$, $R=\sqrt{R_1 R_2}$.  
Given this the  main result of \cite{kpv2018} is the following:
\begin{theorem}\label{theorem liu}
Under the hypothesis and with the notation of Lemma \ref{lemma free} there exist a sequence $\lambda_n\to 0$,  and a sequence of  maximal solutions $u_n$   of the Liouville problem (\ref{liu1}) 
with the following properties:
\begin{itemize}
\item[(i)]
It holds
\begin{align}
\label{liu 5}
\frac{u_{n}}{2\log\frac{1}{\lambda_n}}\longrightarrow H^\pm_\gamma,\quad \mbox{as}\ \lambda_n\to 0,
\end{align}

over compact subsets of $\Omega^\pm$.
\item[(ii)]
We have
\begin{align}
\label{liu 6}
\frac{\lambda_n^2}{2\log\frac{1}{\lambda_n}}\int_\Omega e^{\,u_{n}}\,dx\longrightarrow \frac{4\pi}{\log{\sqrt{\frac{R_1}{R_2}}}}\quad \mbox{as}\ \lambda_n\to 0
\end{align}
\end{itemize}
\end{theorem}
By analogy with the expression in (\ref{mass1}) 
we interpret the right hand side of  (\ref{liu 6}) as the mass of the blow up  of the maximal solution. Note that it depends solely  on the conformal class of $\Omega$.

Our result complements \cite{nagasaki1990}  as it shows that at least in the  case of general doubly connected domains {case (iii) occurs} in general. What's more, we also describe the way that the whole domain blow up actually happens. First, we find the curve of concentration of the blow up in terms of the free boundary problem (\ref{liu 3})--(\ref{liu 4}), second we find the exact form of the line bubble in terms of the one dimensional solution of the Liouville equation (see section \ref{subsec 1} below), third we describe how this local behavior is mediated with the far field approximation of the solution given by the scaled  harmonic measures $H^\pm_\gamma$. 

The problem of determining  the curve $\gamma$ is interesting in its own right. The case of doubly connected domains is relatively simple because of the conformal equivalence between $\Omega$ and an annulus. For general multiply connected domains solving (\ref{liu 3})--(\ref{liu 4}) appears to be more complicated but in fact it can be done quite easily, see Remark \ref{rem 1} below. For brevity here we will  discuss the former case only.

The importance of (\ref{liu1})  
lies in the fact that  there are several important problems in geometry and physics for which the Liouville equation is a simple model problem. One of them, that we will discuss later (in section \ref{mean field})  is the mean field equation. Another related problem we should mention (but will not address here) is the   prescribed Gaussian curvature  equation  (known as the Nirenberg problem in the case of the sphere) -- we refer to  \cite{10.2307/1970993, chang1987,chang2004non} and the references therein.  Recently, a version of this problem, where additionally the  geodesic curvature on the boundary is also prescribed  has been studied in \cite{lsmr2018}.

Our work was in part inspired by \cite{gladiali2007} where (\ref{liu 6}) of our theorem was proven for the annuli. Another closely related results are contained in \cite{dpw2006, pv, DELPINO20163414}  where solutions of the stationary Keller-Segel system concentrating on the boundary of the domain were constructed (see also the discussion in section \ref{sec ks}).

%

\subsection{Existence of the free boundary}\label{sec free}

In this section we will prove Lemma \ref{lemma free}.  Under the assumptions of Lemma \ref{lemma free} it is known that $\Omega$ is conformally equivalent to an annulus (see Theorem 4.2.1 in \cite{krantz_book1}). Thus we have a holomorphic, bijective map $\psi\colon\Omega\to B_{R_1}\setminus B_{R_2}$ with some $R_1>R_2>0$. To show Lemma \ref{lemma free} we first solve the problem of finding the curve $\gamma$ in the annulus $A=B_{R_1}\setminus B_{R_2}$ and then pull it back to $\Omega$ using $\psi$. As a candidate for $\gamma$ we take $C_R=\{|x|=R\}$, where $R$ will be adjusted to satisfy the free boundary problem.  We will denote by $A^+=\{R_1>|x|>R\}$ and $A^-=\{R>|x|>R_2\}$. Functions $H^\pm_{C_R}$ should satisfy the following set of conditions
\begin{equation}\label{an free}
\begin{aligned}
\Delta H^\pm_{C_R}&=0, \quad \mbox{in}\quad  A^\pm,\\
H^\pm_{C_R}&=0, \quad \mbox{on}\quad  \partial A^\pm\cap \partial A,\\
H^\pm_{C_R}&=1, \quad \mbox{on}\quad C_R,\\
\partial_r H^+_{C_R}+\partial_r H^-_{C_R} &=0, \quad \mbox{on}\quad C_R.
\end{aligned}
\end{equation}
It is rather easy to see that we should have
\[
H^\pm_{C_R}=a^\pm+b^\pm \log r,
\]
and that all conditions in (\ref{an free}) will be satisfied when $R=\sqrt{R_1 R_2}$ and
\begin{equation}
\label{abplus}
\begin{aligned}
a^+=-\frac{\log R_1}{\log\left(\sqrt{\frac{R_2}{R_1}}\right)}, &\quad a^-=-\frac{\log R_2}{\log\left(\sqrt{\frac{R_1}{R_2}}\right)},\\
b^+=\frac{1}{\log\left(\sqrt{\frac{R_2}{R_1}}\right)}, &\quad b^-=\frac{1}{\log\left(\sqrt{\frac{R_1}{R_2}}\right)}.
\end{aligned}
\end{equation}
We let $\gamma=\psi^{-1}(C_R)$ and 
\[
H^\pm_{\gamma}=H^\pm_{C_R}\circ \psi.
\]
Since $\psi$ is a conformal map it is evident that $H^\pm_\gamma$ satisfies the assertions of Lemma \ref{lemma free}. Observe that by definition we have
\[
\partial_n H^-_\gamma=|\partial_n\psi|\partial_r H^-_{C_R}\circ \psi=-|\partial_n\psi|\partial_r H^+_{C_R}\circ\psi=-\partial_n H^+_\gamma.
\]
Hence, along $\gamma$ 
\begin{equation}
\label{hplusminus}
\partial_n H^-_\gamma=|\partial_n\psi| \frac{b^-}{\sqrt{R_1 R_2}}, \qquad \partial_n H^+_\gamma=|\partial_n\psi| \frac{b^+}{\sqrt{R_1 R_2}}.
\end{equation}

\subsection{Derivation of the free boundary problem}\label{subsec 1} 

In this section we will consider  a slightly modified version of (\ref{liu1}), namely
\begin{equation}
\label{liu 1}
\begin{aligned}
\Delta u+\lambda^2 V(x)e^{\,u}&=0, \quad\mbox{in}\ \Omega\subset \R^2,\\
u&=0, \quad\mbox{on}\, \partial\Omega,
\end{aligned}
\end{equation}
where the potential $V(x)$  is a smooth, positive  function in $\Omega$. At first sight it may seem that the introduction of $V$ should  alter somehow the conditions (\ref{liu 3})--(\ref{liu 4}) but as we will see this is not the case. The reason is that the effect of the potential is of order $\mathcal O(1)$ while the free boundary problem is seen at the order $\mathcal O(\log\frac{1}{\lambda})$.

We fix a smooth, simple closed curve $\gamma\subset \Omega$, a candidate for the free boundary.
To describe   the Fermi  coordinates of the curve $\gamma$ let us  denote the arc length  parametrization of $\gamma$ by $s$ and let $t(x)=\mathrm{dist}(\gamma, x)$ to be the signed distance to $\gamma$ chosen in agreement with its orientation so that 
\[
x=\gamma(s)+t n(s).
\]
For every point $x$ sufficiently close to $\gamma$,  the map
$x\longmapsto (s, t)$
is a diffeomorphism. We will denote this diffeomorphism by $X_\gamma$, so that $X_\gamma(x)=(s,t)$. In what follows we will often express functions globally defined in 
$\Omega$ or $\Omega^\pm$ in terms of the Fermi coordinates keeping in mind that these expressions are correct only when $|\mathrm{dist}(x,\gamma)|<\delta$ with some $\delta>0$ small. 

To find the approximation of the maximal solution near $\gamma$ let us consider  the following ODE:
\begin{align}
\label{1d profile}
\begin{aligned}
&u''+e^{\,u}=0, \quad \mbox{in}\ \R,\\
& u'(0)=0.
\end{aligned}
\end{align}
This problem has an explicit solution
\begin{align}\label{u}
U(t)=\log\left({2\sech^2 t}\right),
\end{align}
whose asymptotic behavior is given by:
\begin{align}\label{u-asy}
U(t)=-a_0|t|+b_0 +\mathcal{O}(e^{\,-a_0|t|}), \quad |t|\to \infty, \quad\mbox{where}\  a_0=2\ \hbox{and}\ b_0=\log 2.
\end{align}

Considering the equation (\ref{liu1}) 
we observe that if $v$ is a solution of  $\Delta v+e^{\,v}=0$
then $u(x)=v(\lambda\mu x)+2\log\mu$ is, for any constant  $\mu>0$, a solution of $\Delta u +\lambda^2 e^{\,u}=0$. 
We introduce the scaling function
\begin{align*}
\mu_\lambda\colon \gamma\to \R_+,
\end{align*}
to be determined later on.  {\it A priori} we assume only
\begin{align*}
\mu_\lambda=\mathcal{O}\left(\frac{\log\frac{1}{\lambda}}{\lambda}\right).
\end{align*}
We also need to take into account the potential so we set $h_\gamma(s)=V\circ X^{-1}_\gamma(s, 0)$ and define the approximate solution by
\begin{equation}
\label{uzero 2}
v_0(x)=U\left(\lambda\mu_\lambda t\right)+2\log\mu_\lambda-\log h_\gamma, \quad x=X_{\gamma}^{-1}(s,t).
\end{equation}
Note that $v_0$ is  well defined as a function of 
$x\in\Omega\cap\{|\mathrm{dist}\, (x,\gamma)|<\delta\}$ with some $\delta>0$. 
By \eqref{u} and \eqref{u-asy} we deduce that 
\begin{align}
\label{uzero 3}
v_0 \circ X^{-1}_\gamma (s,t)=-a_0 \lambda\mu_\lambda |t|+b_0+2\log\mu_\lambda-\log h_\gamma+\mathcal O\left(e^{-a_0\lambda\mu_\lambda|t|}\right).
\end{align}

Next we consider  the outer approximation.   In each of the components $\Omega^\pm$ of $\Omega\setminus \gamma$ we define the outer  approximation $w_0^\pm$ by:
\begin{align}
\label{ext 1}
\begin{aligned}
\Delta w_0^\pm&=0, \quad\mbox{in}\ \Omega^\pm, \\
w_0^\pm&=0, \quad \mbox{on}\ \partial\Omega^\pm\cap\partial \Omega.
\end{aligned}
\end{align}
Note that in this problem we are missing a boundary condition on $\gamma$ for each of the unknown functions. This boundary condition will be determined through a matching condition between the inner approximation $v_0$ and the outer approximation $w_0^\pm$ and it will provide eventually the free boundary condition (\ref{liu 4}). To proceed let us express the outer approximations in terms  of the Fermi coordinates of $\gamma$, and then expand them formally in terms of the variable $t$:
\begin{align*}
w_0^\pm\circ X^{-1}_\gamma (s,t)&=w_0^\pm \circ X^{-1}_\gamma (s, 0)+t\partial_t w_0^\pm \circ X^{-1}_\gamma(s,t)\mid_{t=0}+\dots\\
&=w_0^\pm (s, 0)+t \partial_n w_0^\pm \circ X^{-1}_\gamma(s,0)+\dots
\end{align*}
Next we let $\eta$ to be the {\it inner  variable}
\begin{equation}
\label{def eta 1}
\eta=\lambda\mu_\lambda t \Longrightarrow t=\frac{\eta}{\lambda\mu_\lambda}
\end{equation}
which is more convenient to express the matching conditions.   
To match the inner and outer approximations  the following identities should  hold:
\begin{align}
\label{match 1}
\begin{aligned}
 w_0^+\circ X^{-1}_\gamma(s,0)+\frac{\eta}{\lambda\mu_\lambda}\left(\partial_n w_0^+\circ X^{-1}_\gamma\right)(s,0)&=-a_0\eta+b_0+2\log\mu_\lambda-\log h_\gamma,\\
w_0^- \circ X^{-1}_\gamma(s,0)+\frac{\eta}{\lambda\mu_\lambda}\left(\partial_n w_0^-\circ X^{-1}_\gamma\right)(s,0)&=a_0\eta+b_0+2\log\mu_\lambda-\log h_\gamma.
\end{aligned}
\end{align}
This leads to the following conditions on $\gamma$
\begin{equation}
\label{match 1a}
 \begin{cases}
w_0^+=b_0+2\log\mu_\lambda-\log h_\gamma,\\
\partial_n w_0^+ =-a_0\lambda\mu_\lambda,
\end{cases} \ \hbox{and}\quad
 \begin{cases}
w_0^-=b_0+2\log\mu_\lambda-\log h_\gamma,\\
\partial_n w_0^-=a_0\lambda\mu_\lambda.
\end{cases}
 \end{equation}
These boundary conditions together with the boundary conditions on $\partial \Omega\cap \partial \Omega^\pm$ give an overdetermined, nonlinear  problem for $\mu_\lambda$, which seems to be rather difficult. To avoid this complication we note that  at this point we only need to satisfy  conditions (\ref{match 1a})  with certain precision. For now it suffices that the difference between the left and the right hand sides is of order $\mathcal O\left(\frac{\log\log\frac{1}{\lambda}}{\log \frac{1}{\lambda}}\right)$ in the matching of $w_0^\pm$ and $\mathcal O(1)$ in the matching of $\partial_n w_0^\pm$.  To accomplish this we introduce the following  scaling function of the small parameter $\lambda$:
\begin{align}
\beta &=2\log\frac{1}{a_0\lambda}+b_0,\quad a_0=2, \quad b_0=\log 2, \label{beta}\\
\mu_\lambda&=-\frac{(\beta+2\log\beta)}{a_0\lambda}\partial_n H_\gamma^+=\frac{(\beta+2\log\beta)}{a_0\lambda} \partial_n H_\gamma^- >0. \label{mu}
\end{align}
With these definitions we set
\begin{equation}\label{wzero}
w_0^\pm=(\beta+2\log\beta) H_\gamma^\pm+\tilde H^\pm,
\end{equation}
where
\[
\begin{aligned}
\Delta \tilde H^\pm&=0, \qquad \mbox{in}\ \Omega^\pm,\\
\tilde H^\pm &=0, \qquad \mbox{on}\ \partial\Omega^\pm\cap\partial\Omega,\\
\tilde H^\pm&=2\log |\partial_nH_\gamma^\pm|-\log h_\gamma, \qquad \mbox{on}\ \gamma.
\end{aligned}
\]
We check 
\begin{equation}
\label{match 1b}
\partial_n w_0^\pm\pm a_0\lambda\mu_\lambda=\partial_n \tilde H^\pm=\mathcal O(1)\ \hbox{on}\ \gamma,
\end{equation}
so that the matching conditions for the derivatives are satisfied as we wanted. We claim that  
\begin{equation}
\label{match 1c}
 w_0^\pm-b_0-2\log\mu_\lambda+\log h_\gamma=\mathcal O\left(\frac{\log\log\frac{1}{\lambda}}{\log\frac{1}{\lambda}}\right)\ \hbox{on}\ \gamma 
\end{equation}
as needed.
Indeed
\[
\begin{aligned}
w_0^+-b_0-2\log\mu_\lambda+\log h_\gamma&=(\beta+2\log \beta)H^+_\gamma+2\log(-\partial_n H^+_\gamma)-\log h_\gamma-b_0-2\log\left(\frac{1}{a_0\lambda}\right)
\\
&\quad-2\log\beta-2\log(-\partial_n H^+_\gamma)+\log h_\gamma-2\log\left(1+\frac{2\log\beta}{\beta}\right)\\
&=-2\log\left(1+\frac{2\log\beta}{\beta}\right)\\
&= \mathcal O\left(\frac{\log\beta}{\beta}\right)=\mathcal O\left(\frac{\log\log\frac{1}{\lambda}}{\log\frac{1}{\lambda}}\right) ,
\end{aligned}
\]
(we have  used $H_\gamma^+=1$ on $\gamma$). Of course we compute similarly the error of $w_0^--b_0-2\log\mu_\lambda$. 

The asymptotic formula for the maximal solution, together with the explicit form of the functions $H^\pm_\gamma$ can be used to calculate (formally)  the mass of the maximal solution as in Theorem \ref{theorem liu} (ii). To begin consider the exponential of the solution $u$. Given $\delta>0$ small we have
\[
\lambda^2e^{\,u}\approx\begin{cases} \lambda^2\exp\left(U\left(\lambda\mu_\lambda t\right)+2\log\mu_\lambda-\log h_\gamma\right), \quad |t|<\delta,\medskip \\
\lambda^2\exp\left((\beta+2\log\beta) H_\gamma^\pm+\tilde H^\pm\right), \quad \mathop{dist}(x, \gamma)>\delta.
\end{cases}
\]  
It can be shown that there exists a constant $c_0>0$  such that  we have 
\[
H_\gamma^\pm(x)\leq 1-c\delta, \quad \mathop{dist}(x, \gamma)\geq \delta.
\]
It follows that, when $\delta = \delta_\lambda=\frac{M\log\beta}{\beta}$ with some $M>0$ large, then
\[
\lambda^2\exp\left((\beta+2\log\beta) H_\gamma^\pm+\tilde H^\pm\right)\lesssim e^{-c_0M\log\beta}=o(\beta),\quad \mathop{dist}(x, \gamma)>\delta_\lambda.
\]
Using this we get
\[
\frac{\lambda^2}{\beta} \int_\Omega V(x)e^{\,u}\,dx=\frac{\lambda^2}{\beta}\int_{|t|<\delta_\lambda}V(x)e^{\,u}\,dx +o(1),
\]
or in other words the mass of the maximal solution is concentrated near $\gamma$. Expressing the integral above in the local coordinates of $\gamma$ we get
\begin{equation}
\label{mass 1}
\begin{aligned}
\frac{\lambda^2}{\beta}\int_{|t|<\delta_\lambda}V(x)e^{\,u}\,dx&\approx \frac{\lambda^2}{\beta}\int_{|t|<\delta_\lambda} h_\gamma(s) \exp\left(U\left(\lambda\mu_\lambda t\right)+2\log\mu_\lambda-\log h_\gamma\right)\,dt ds\\
&\approx \frac{\lambda^2}{\beta} \int_0^{|\gamma|}\!\int_{-\delta_\lambda}^{\delta_\lambda} \mu_\lambda^2\exp\left(U(\lambda\mu_\lambda t)\right)\,dt ds\\
&\approx \int_0^{|\gamma|}\!\int_{-\infty}^{\infty}\frac{\lambda\mu_\lambda}{\beta} e^{\,U(\eta)}\,d\eta ds\\
&\approx \int_0^{|\gamma|} \frac{1}{2} |\partial_n H^\pm_\gamma(s,0)|\,ds\int_{-\infty}^\infty \exp\left(\log(2\sech^2 \eta)\right)\,d\eta \\
&\approx \frac{2|b^-|}{\sqrt{R_1 R_2}}\int_0^{|\gamma|} |\partial_n\psi| \,ds\\
&=\frac{4\pi R}{\sqrt{R_1 R_2}\log\sqrt{\frac{R_1}{R_2}}}\\
&=\frac{4\pi}{\log{\sqrt{\frac{R_1}{R_2}}}}.
\end{aligned}
\end{equation}
Above we have used  that $|\partial_n\psi|=|\partial_\tau\psi|$ by the Cauchy-Riemann equations and 
\begin{equation}
\label{sech identity}
\int_{-\infty}^\infty \exp\left(\log(2\sech^2 \eta)\right)\,d\eta=2\int_{-\infty}^\infty \sech^2 \eta\,d\eta=4.
\end{equation}

\begin{remark}
\label{rem 1}
It is evident that the free boundary problem (\ref{liu 3})--(\ref{liu 4})  can be  stated when $\Omega$ is an arbitrary bounded and smooth domain in $\R^2$. Suppose that  we can find $W^\pm_\gamma$ in such general case and let us define 
\[
u_\gamma=\begin{cases} W^+_\gamma, \quad \mbox{in} \quad \Omega^+,\\
2-W^-_\gamma, \quad \mbox{in} \quad \Omega^-.
\end{cases}
\]
Then in fact $u_\gamma$ is a harmonic function in the whole $\Omega$ such that $u_\gamma=0$ on this  part of $\partial\Omega$ where $W^+_\gamma=0$  and $u_\gamma=2$ 
on this part of $\partial\Omega$ where $W^-_\gamma=0$. At the same time $u_\gamma=1$ on $\gamma$. Conversely, given two disjoint components $\partial\Omega_\ell$, 
$\ell=0,2$ such that $\partial\Omega=\partial\Omega_0\cup\partial\Omega_2$ we can find a solution of the free boundary problem  (\ref{liu 3})--(\ref{liu 4}) by  solving
\begin{equation}
\label{fh 1}
\begin{aligned}
\Delta u&=0, \quad\mbox{in} \quad \Omega,\\
u&=\ell,  \quad\mbox{in} \quad \partial\Omega_\ell,\quad  \ell=0,2,
\end{aligned}
\end{equation}
and defining $\gamma=\{u=1\}$.  This way the problem (\ref{fh 1}) leads to a kind of generalized solution of the free boundary problem. Note that we can  also find $\gamma$ by minimizing $\int_\Omega |\nabla u|^2$ among all $H^1(\Omega)$  functions such that $u=\ell$, on $\partial\Omega_\ell$, $\ell =0,2$. With this definition several solutions of the free boundary problem can be possibly found for domains with sufficiently rich topology. 
\end{remark}

\section{The mean field equation}\label{mean field}
\setcounter{equation}{0}
\subsection{Derivation of the free boundary problem}
Liouville's equation (\ref{liu 1}) belongs to a larger family of problems, one of them is the mean field model
\begin{equation}
\label{mean 1}
\Delta_g u+\rho\left (\frac{V(x) e^{\,u}}{\int_M V e^{\,u}\,d\mu}-c_M\right)=0,\qquad c_M=\frac{1}{\mathop{vol}(M)},
\end{equation}
on a compact, two dimensional, closed  Riemannian manifold $(M,g)$ assuming $V>0$. This equation appears in statistical mechanics \cite{Caglioti1992, Caglioti1995, doi:10.1002/cpa.3160460103} and Chern-Simmons-Higgs theory \cite{PhysRevLett.64.2234,PhysRevLett.64.2230} and its  most complete existence theory was developed by Chen and Lin \cite{doi:10.1002/cpa.3014,doi:10.1002/cpa.10107}, see also \cite{Li1999,DING1999653,Struwe1998, MR2483132,MR2671126, MR2409366,MR2456884,MR3263503}. In very general terms, given suitable  assumptions on $M$,  these results show the existence of bubbling solutions as the parameter $\rho\to \infty$ with the number of bubbles increasing to infinity as 
$\rho$ crosses the values $\rho_m=8\pi m$, $m\in \mathbb N$.  For such solutions we have 
\[
\int_M V e^{\,u}\,d\mu\longrightarrow 8\pi m, \qquad \rho\nearrow 8\pi m.
\]
Note that the mass of  the bubbling solutions increases as $\rho$ as $\rho\to \infty$. 
Here we will consider the maximal solutions meaning that 
\begin{equation}
\label{mass increase}
\frac{\rho}{\int_M V e^{\,u}\,d\mu}=\lambda^2\longrightarrow 0, \qquad \rho\to \infty.
\end{equation}
In contrast with the bubbling solutions the mass of the maximal solutions increases faster than $\rho$ as $\rho \to \infty$.
Calculations of the previous section suggest that we could take $\rho=\mathcal{O}(\log\frac{1}{\lambda})$ but strictly speaking it is not necessary so that we will assume only that 
\[
\rho=\rho_\lambda\longrightarrow +\infty, \quad \lambda\searrow 0,
\]
and (\ref{mass increase}) holds. 
With this in mind (\ref{mean 1}) becomes
\begin{equation}
\label{mean 2}
\Delta_g u+\lambda^2 V(x) e^{\,u}=\rho_\lambda c_M. 
\end{equation} 
To simplify the discussion we will  consider (\ref{mean 2}) in a bounded, smooth  domain $\Omega\subset \R^2$,  assuming homogeneous Neumann boundary conditions and $V\equiv 1$:
\begin{equation}
\begin{aligned}
\Delta u+\lambda^2 e^{\,u}&=\rho_\lambda c_\Omega, \quad c_\Omega=\frac{1}{|\Omega|},\\
\partial_n u&=0, \quad \mbox{on}\quad  \partial \Omega.
\end{aligned}
\label{mean 3}
\end{equation}
Integrating (\ref{mean 3}) and taking into account the boundary conditions we see that we should have
\[
\frac{\lambda^2}{\rho_\lambda}\int_\Omega V(x) e^{\,u}=1,
\]  
which shows that (\ref{mean 3}) is a reasonable model problem for (\ref{mean 1}). 

We proceed to define the approximate solution for (\ref{mean 3}).  As before, we choose a smooth, closed  curve $\gamma$  in $\Omega$ and define the inner approximation by (\ref{uzero 2}). At  this moment $\mu_\lambda$ is still unknown.  At the distance $\mathcal O(1)$ from $\gamma$ we expect that 
\[
u\sim \rho_\lambda, \qquad \lambda^2 e^{\,u}=o(\rho_\lambda)
\]
hence we should have $u\approx w_0^\pm$, where $w^\pm_0= \rho_\lambda c_\Omega K^\pm_\gamma$ and 
\begin{equation}
\label{mean 4}
\begin{aligned}
\Delta K_\gamma^\pm&=1, \quad \mbox{in}\ \quad\Omega^\pm,\\ 
\partial_n K_\gamma&=0\quad \mbox{on}\ \quad\partial \Omega^\pm\cap \partial\Omega.
\end{aligned}
\end{equation}
To find  the free boundary conditions we use exactly the same argument of matching the inner and the outer approximations as in the previous case, see the relations in (\ref{match 1a}).  As a consequence  we need 
\begin{equation}
\begin{aligned}
\label{mean 5}
K^\pm_\gamma&=1, \quad \mbox{on}\quad \gamma,\\
\partial_n K^+_\gamma+\partial_n K^-_\gamma&=0, \quad \mbox{on}\quad \gamma.
\end{aligned}
\end{equation}
Keep in mind that $n$ is the unit normal exterior to $\Omega^-$, by the choice of the orientation on $\gamma$. Suppose that (\ref{mean 4})--(\ref{mean 5}) can be solved (we will provide an example in section \ref{sec disc}). We need
\[
\rho_\lambda c_\Omega \partial_n K^-_\gamma=a_0\lambda\mu_\lambda=-\rho_\lambda c_\Omega \partial_n K^+_\gamma, \quad a_0=2,
\]
hence the definition of $\mu_\lambda$, c.f (\ref{mu}):
\begin{equation}
\label{mean 6}
\mu_\lambda =\frac{\rho_\lambda c_\Omega |\partial_n K^\pm_\gamma|}{a_0\lambda}.
\end{equation}
Using this approximation we calculate the mass of the solution $\lambda^2 \int_\Omega e^{\,u}\,dx$ as in (\ref{mass 1}). By  equation (\ref{mean 3}) we have of course
\begin{equation}
\label{mass 2}
\frac{\lambda^2}{\rho_\lambda} \int_\Omega e^{\,u}\,dx =1.
\end{equation}
On the other hand following (\ref{mass 1}) (calculation is essentially the same up to the fourth line replacing $\beta$ by $\rho_\lambda$ and $H^\pm_\gamma$ by $K^\pm_\gamma$)
\begin{equation}
\label{mass 3}
\begin{aligned}
\frac{\lambda^2}{\rho_\lambda} \int_\Omega e^{\,u}\,dx &\approx \int_0^{|\gamma|} \frac{1}{2} |\partial_n K^\pm_\gamma(s,0)|\,ds\int_{-\infty}^\infty \exp\left(\log(2\sech^2 \eta)\right)\,d\eta \\
&=2c_\Omega\int_\gamma |\partial_n K^\pm_\gamma|=2c_\Omega |\Omega^\pm|=1,
\end{aligned}
\end{equation}
where we have used (\ref{sech identity}) and
\[
|\Omega^+|=\int_{\Omega^+}\Delta K^+=-\int_\gamma \partial_n K^+=\int_\gamma \partial_n K^-=\int_{\Omega^+}\Delta K^+=|\Omega^-|,
\]
and $|\Omega^+|+|\Omega^-|=|\Omega|$. Thus our approximation is consistent with (\ref{mass 2}). 

Let us discuss  existence  of  another set where the mass could concentrate.  One possibility  is suggested by \cite{DELPINO20163414}: the mass is concentrated on $\partial \Omega$. To take this into account suppose that $\alpha=\partial\Omega$ is a simple, smooth closed curve. Let us suppose that $\alpha$ is contained in the same connected component of $\R^2\setminus \gamma$ as $\Omega^-$. Let now  the function $K_{\gamma}^-$ be the solution of 
\begin{equation}
\label{mean 7}
\begin{aligned}
\Delta K^-_{\gamma}&=1,\quad \mbox{in}\quad  \Omega,\\
K^-_\gamma &=1, \quad \mbox{on}\quad  \partial \Omega=\alpha,\\
K^-_\gamma &= 1, \quad \mbox{on}\quad \gamma.
\end{aligned}
\end{equation}
The function $K^+_\gamma$ is the  solution of (\ref{mean 4}) as before. 

We associate to $\alpha$ local coordinates $(\sigma, \tau)$ where $\sigma$ is the arc length and $\tau>0$ the distance. Consider a scaling function $\nu_\lambda$ as before and let the boundary approximation of the maximal solution to be
\[
U(\lambda \nu_\lambda\tau)+2\log\nu_\lambda,
\] 
and the inner approximation of the maximal solution to be
\[
U(\lambda \mu_\lambda t)+2\log\mu_\lambda.
\]
As before, we keep in mind that both scaling functions $\mu_\lambda, \nu_\lambda$ are yet to be determined. 
Next, let the outer approximation be defined at the distance larger than $O(\delta_\lambda)$ where $\delta_\lambda=\frac{M\log\rho_\lambda}{\rho_\lambda}$ from $\alpha$ or $\gamma$ by 
\[
w_0= \begin{cases} w_0^-=\rho_\lambda c_\Omega K_\gamma^-, \quad \mbox{in}\ \quad \Omega^-,  \medskip \\
w_0^+= \rho_\lambda c_\Omega K^+_\gamma, \quad\mbox{in}\quad \Omega^+.\\
\end{cases}
\]

To determine $\nu_\lambda$ we use the matching condition of the form (\ref{match 1a}) for the derivative:
\begin{equation}
\label{mean 8}
\partial_n w^-_0=a_0\lambda\nu_\lambda\Longrightarrow \nu_\lambda=\frac{\rho_\lambda c_\Omega \partial_n K_\gamma^-}{a_0\lambda}. 
\end{equation}
To find $\mu_\lambda$ we argue similarly. On $\gamma$ we need:
\begin{equation}
\label{mean 9}
\begin{aligned}
\partial_n w_0^+ &=-a_0\lambda\mu_\lambda,\\
\partial_n w_0^- &=a_0\lambda\mu_\lambda,
\end{aligned}
\end{equation}
hence in addition to $K_\gamma^\pm=1$ on $\gamma$  we should require
\begin{equation}
\label{mean 10}
\partial_n K^+_\gamma +\partial_n K^-_\gamma=0.
\end{equation}
Finally, we get
\begin{equation}
\label{mean 11}
\mu_\lambda=\frac{\rho_\lambda c_\Omega |\partial_n K^\pm_\gamma|}{a_0\lambda}
\end{equation}
Now, we can calculate the mass of the maximal solution by similar calculations as those in (\ref{mass 1}) and (\ref{mass 3}). Note that  the integral over $\alpha$ will carry a factor of $\frac{1}{2}$ since in the volume integral we only integrate on one side of $\alpha$. The result is
\begin{equation}
\label{mass 4}
\begin{aligned}
\frac{\lambda^2}{\rho_\lambda} \int_\Omega e^{\,u}\,dx &=2c_\Omega\int_\gamma |\partial_n K^\pm_\gamma|+c_\Omega\int_\alpha \partial_n K_\alpha\\
&= c_\Omega |\Omega^+|+c_\Omega|\Omega^-|\\
&=1,
\end{aligned}
\end{equation}
as it should be.

\subsection{Examples}
\subsubsection{Problem (\ref{mean 3}) in the disc}\label{sec disc}

When $\Omega=B_{R_1}$   the solution of the free boundary problem (\ref{mean 4})--(\ref{mean 5}) should be  $\gamma=\{|x|=R\}$ with $R=R_1/\sqrt{2}$, since it divides $\Omega$ into two regions of equal areas.   We will verify this solving directly (\ref{mean 4})--(\ref{mean 5}). To this end we suppose that $K^\pm_\gamma$ are radial functions so that  
\[
K^+_\gamma(r)= a^++\frac{1}{4}r^2, \qquad K^-_\gamma(r)=a^-+b^-\log r+ \frac{1}{4}r^2.
\]
The boundary conditions and the free boundary conditions give the equations
\begin{equation}
\label{ex free 1}
\begin{aligned}
a^++\frac{1}{4}R^2&=1,\\
a^-+b^-\log R+\frac{1}{4}R^2&=1,\\
b^-+\frac{1}{2}R_1^2&=0,\\
R^2+b^-=0.
\end{aligned}
\end{equation}
For the third and the fourth equations to be consistent we need
\[
R=\frac{R_1}{\sqrt{2}},
\]
as expected.

Next, consider the modification of the free boundary problem as in (\ref{mean 7}). In this case we no longer have $|\Omega^+|=|\Omega^-|$ and so it is not so easy to guess what the free boundary $\gamma$ should be. Modifying (\ref{ex free 1}) to account for the new boundary condition on $K^-_\gamma$ we get 
(compared with (\ref{ex free 1}) only the third equality is different)
\[
\begin{aligned}
a^++\frac{1}{4}R^2&=1\\
a^-+b^-\log R+\frac{1}{4}R^2&=1,\\
a^-+b^-\log R_1+\frac{1}{4}R_1^2&=1, \\
R^2+b^-&=0.
\end{aligned}
\]
Eliminating $a^-,b^-$ in the last three equations and  setting $r=\frac{R^2_1}{R^2}$ we see that the problem is to find an $r\in (1,\infty)$ such that 
\[
\frac{1}{2}\log r-\frac{1}{4} r=-\frac{1}{4}.
\]
For the function $f(r)=\frac{1}{2}\log r-\frac{1}{4} r$ we have 
\[
f(1)=-\frac{1}{4}, \quad f'(r)>0, \quad r\in (1,2), \quad f'(r)<0, \quad r>2
\]
hence there exists $r^*>2$ such that $f(r^*)=-\frac{1}{4}$ and thus
\[
R=\frac{R_1}{\sqrt{r^*}}
\]
gives the radius of $\gamma$ in this case.
Note however that these two examples  are not the only possible solutions for we could take as $\gamma$ the diameter, find $K^-_\gamma$ in the half-disc and then define $K^+_\gamma$ by the even reflection.  This suggests that there is a multitude of the  maximal solutions  in general. 

\subsubsection{The free boundary problem for the  mean field model on a manifold}
Based on the derivation of the free boundary problem for the Liouville equation and for the model problem (\ref{mean 3}) it is easy to guess what the free boundary problem associated with the maximal solutions of (\ref{mean 1}) should be. To state it in precise and more general terms let $\Omega^\pm\subset M$ be two open sets such that 
$\overline{\Omega^+\cup\Omega^-}=M$. The free boundary is the common boundary of the two sets $\gamma=\partial\Omega^\pm$. The problem of  determining $\gamma$ amounts to finding  $K^\pm_\gamma$ such that 
\begin{equation}
\label{free bdry mf}
\begin{aligned}
\Delta_g K^\pm_\gamma&=1, \quad \mbox{in}\quad \Omega^\pm,\\
K^\pm_\gamma&=1, \quad \mbox{on}\quad \gamma,\\
\partial_n K^-_\gamma+\partial_n K^+_\gamma&=0, \quad \mbox{on}\quad \gamma,
\end{aligned}
\end{equation} 
where $n$ is the unit outer normal vector on $\Omega^-$.  
Problem (\ref{free bdry mf}) does not seem to be as simple as the problem considered in (\ref{liu 3})--(\ref{liu 4}). Note that now using conformal map to transfer this problem to some standard domain (e.g. the sphere) does not work.  On the other hand it is easy to give examples of solutions for some concrete manifolds.  For example if $M=S^2$ then $\Omega^\pm$ should be the half spheres. If $M={\mathbb T}^2$ (two dimensional torus) then two possible solutions come  come to mind immediately. Agreeing that the axis of the rotation of the torus is the $z$ axis one solution is to take $\Omega^+= {\mathbb T}^2\cap \{x>0\}$, $\Omega^-= {\mathbb T}^2\cap \{x<0\}$. Another is to take $\Omega^+={\mathbb T}^2\cap \{z>0\}$, $\Omega^-={\mathbb T}^2\cap \{z<0\}$. Again we see that the solution to the free boundary problem  is not unique, even if we mod out the obvious symmetry. 

\subsubsection{Stationary solutions of the Keller-Segel problem}\label{sec ks}
The stationary Keller-Segel on a bounded domain $\Omega \subset \R^2$  can be written in  the following form
\begin{equation}
\label{ks 1}
\begin{aligned}
\Delta u+u&=\lambda e^{\,u}, \quad \mbox{in}\quad \Omega,\\
\partial_n u&=0, \quad \mbox{on}\quad \partial \Omega.
\end{aligned}
\end{equation}
where $\lambda>0$ is a small parameter (see for instance \cite{dpw2006}). 
Solutions blowing up in the whole domain $\Omega$ and with the mass concentrating on the boundary were constructed in \cite{DELPINO20163414}. In the case of radially symmetric domains ($\Omega$ is a disc or an annulus) the analogs of maximal solutions concentrating on internal circles of $\Omega$ were found by bifurcation analysis combined with  variational methods in \cite{Bonheure2017, MR3625083,2017arXiv170910471B} (see also \cite{BONHEURE2016455} for a related higher dimensional problem). Formally it is not hard to find the  free boundary problem that would  determine the curve $\gamma\subset \Omega$  of mass concentration in a general domain. It should consist of finding functions $U^\pm_\gamma$ such that 
\begin{equation}
\label{ks 2}
\begin{aligned}
\Delta U^\pm_\gamma+U^\pm_\gamma&=0,\quad\mbox{in}\quad \Omega^\pm,\\
\partial_n U^\pm_\gamma&=0, \quad \mbox{on}\quad \partial\Omega,\\
U^\pm_\gamma&=1, \quad \mbox{on}\quad \gamma.
\end{aligned}
\end{equation}
Indeed this problem in the radially symmetric setting has been stated already in \cite{BONHEURE2016455}. It is important to notice that (\ref{ks 2}) can be also stated when $\gamma$ instead of being a single curve has several components. For example in \cite{MR3625083, 2017arXiv170910471B} solutions of (\ref{ks 1}) with mass concentrating on concentric circles were found. These circles can be obtained by solving suitably modified (\ref{ks 2}).  Finally we observe that the free boundary problem in this case can be formulated as a partition problem: find $\gamma$ such that $\|u\|_{H^1(\Omega)}^2$ is minimized in the class of functions $u=1$ on $\gamma$. 


\providecommand{\bysame}{\leavevmode\hbox to3em{\hrulefill}\thinspace}
\providecommand{\MR}{\relax\ifhmode\unskip\space\fi MR }
\providecommand{\MRhref}[2]{%
  \href{http://www.ams.org/mathscinet-getitem?mr=#1}{#2}
}
\providecommand{\href}[2]{#2}

\end{document}